

%
%
\input amsppt.sty
\magnification=\magstep1
\hsize = 6.5 truein
\vsize = 9 truein
\NoBlackBoxes
\TagsAsMath
\NoRunningHeads

\catcode`\@=11
\redefine\logo@{}
\catcode`\@=13

\newskip\sectionskipamount
\sectionskipamount = 24pt plus 8pt minus 8pt
\def\sectionskip{\vskip\sectionskipamount}
\define\sectionbreak{%
	\par  \ifdim\lastskip<\sectionskipamount
	\removelastskip  \penalty-2000  \sectionskip  \fi}
\define\section#1{%
	\sectionbreak	
	\subheading{#1}%
	\bigskip
	}

\addto\tenpoint{\normalbaselineskip 18truept \normalbaselines}
\def\label#1{\par%
	\hangafter 1%
	\hangindent 1.25 in%
	\noindent%
	\hbox to 1.25 in{#1\hfill}%
	\ignorespaces%
	}

\define\op#1{\operatorname{\fam=0\tenrm{#1}}} 

\redefine\qed{{\unskip\nobreak\hfil\penalty50\hskip2em\vadjust{}\nobreak\hfil
    $\square$\parfillskip=0pt\finalhyphendemerits=0\par}}

	\define	\x   {\times}
	\let	\< =     \langle
	\let	\> =     \rangle
	\define	 	\a		{\alpha}
	\redefine	\b		{\beta}
	\redefine	\d		{\delta}
	\redefine	\D		{\Delta}
	\define	 	\e		{\varepsilon}
	\define	 	\g		{\gamma}
	\define	 	\G		{\Gamma}
	\redefine	\l		{\lambda}
	\redefine	\L		{\Lambda}
	\define		 \n		{\nabla}
	\redefine	\var	  {\varphi}
	\define		 \s		   {\sigma}
	\redefine	\Sig   {\Sigma}
	\redefine	\t		 {\tau}
	\define	 	\Th		{\Theta}
	\redefine	\O		 {\Omega}
	\redefine	\o	 	{\omega}
	\define		 \z	 	{\zeta}
	\redefine	\i	 	{\infty}
	\define		 \p		 {\partial}

\document
\baselineskip=18 pt

\centerline{\bf Volumes of Restricted Minkowski Sums}
\centerline{\bf and the Free Analogue of the Entropy Power
Inequality}

\vskip .33 truein
\centerline{\smc Stanislaw J. Szarek${}^*$}
\centerline{Department of Mathematics}
\centerline{Case Western Reserve University}
\centerline{Cleveland, Ohio 44106-7058}

\bigskip
\centerline{\smc Dan Voiculescu\footnote"*"{Research supported
  in part by grants from the National Science Foundation}}
\centerline{Department of Mathematics}
\centerline{University of California}
\centerline{Berkeley, California 94720-3840}

\newpage
In noncommutative probability theory independence can be based on
free products instead of tensor products. This yields a highly
noncommutative theory: free probability (for an introduction see
[9]). The analogue of entropy in the free context was introduced by
the second named author in [8]. Here we show that Shannon's entropy
power inequality ([6],[1]) has an analogue for the free entropy
$\chi(X)$ (Theorem 2.1).

The free entropy, consistently with Boltzmann's formula
$S=k\log W$, was defined via volumes of matricial microstates.
Proving the free entropy power inequality naturally becomes a
geometric question.

Restricting the Minkowski sum of two sets means to specify the set
of pairs of points which will  be added. The relevant inequality,
which holds when the set of addable points is sufficiently large,
differs from the Brunn-Minkowski inequality by having the exponent
$1/n$ replaced by $2/n$. Its proof uses the rearrangement
inequality of Brascamp-Lieb-L\"uttinger ([2]). Besides the free
entropy power inequality, note that the inequality for restricted
Minkowski sums may also  underlie the classical Shannon entropy
power inequality (see 3.2 below).

\bigskip
{\bf Acknowledgments}. Part of this work was done while the first
named author was visiting Universit\'e Paris VI and Universit\'e
Marne la Vallee in May--July 1995.  The second named author worked
on this paper while visiting MIT and IHES in 1995.  The authors
express their gratitude to these institutions for their
hospitality and support.

\bigskip\bigskip
{\bf 1. \ The inequality for restricted Minkowski sums}.
If $A,B\subset\Bbb R^n$ (or any vector space), the Minkowski sum of
$A$ and $B$ is defined by
$$
A+B=\{ x+y:(x,y)\in A\times B\} \ .
$$
An important property of the Minkowski sum in $\Bbb R^n$ is the
Brunn-Minkowski inequality\newline ([4],[5])
$$
\l(A+B)^{1/n} \geq \l(A)^{1/n} +\l(B)^{1/n}
$$
where $\l$ denotes $n$-dimensional Lebesgue measures. We introduce
a modified concept of a sum.

\bigskip
{\bf 1.1 \ Definition}. Let $A,B$ be subsets of a vector space and
$\Th\subset A\times B$. We will call
$$
A+_\Th B= \{ x+y:(x,y)\in\Th\}
$$
the restricted (to $\Th$) sum of $A$ and $B$.

\bigskip
We then have the following inequality (in what follows,
all sets and functions are assumed to be measurable; $\l$ denotes
the Lebesgue measure in the appropriate dimension that may vary
from place to place).

\proclaim{1.2 \ Theorem} Let $\rho\in (0,1)$, \ $n\in\Bbb N$
and let $A,B\subset\Bbb R^n$ be such that
$$
\rho\leq\left(\frac{\l(B)}{\l(A)}\right)^{\frac 1n}\leq\rho^{-1} \ .
$$
Furthermore, let $\Th\subset A\times B \subset \Bbb R^{2n}$ be such that
$$
\l(\Th)\geq (1-c\min\{\rho\sqrt{n},1\})\l(A)\l(B) \ .
$$
Then
$$
\l(A+_\Th B)^{2/n} \geq \l(A)^{2/n}+\l(B)^{2/n} \ .  \tag 1.1
$$
($c > 0$ is a numerical constant, independent of $\e$,
$n$, $A$, $B$ and $\Th$.)
\endproclaim

The following simple but illuminating example shows that, in
general, one cannot expect a significantly stronger assertion:
let $B^n$ be the Euclidean ball in $\Bbb R^n$, \
$A=B^n$, $B=\rho B^n$ and $\Th=\{(x,y)\in A\times B:\< x,y\>\leq
0\}$. Then
\roster
\item \  $\l(\Th)=\frac 12\l(A)\l(B)$
\item \  $A+_\Th B=(1+\rho^2)^{\frac 12} B^n$
\endroster
and we have equality in (1.1). We now state a lemma which is an
elaboration of this example

\proclaim{1.3 \ Lemma} Let $\rho,n$ be as in Theorem 1.2 and let
$$
\Th=\{ (x,y):x,y\in\Bbb R^n, |x|\leq 1, \ |y|\leq\rho, \
  |x+y|\leq(1+\rho^2)^{\frac 12}\}
$$
Then
$$
\l(\Th)\leq (1-c\min\{\rho\sqrt{n},1\})\l(B^n)\l(\rho B^n) \ ,
$$
where $c > 0$ is a universal constant.
\endproclaim

We postpone the proof of the lemma (which depends on a careful,
but completely elementary computation) and show how it implies the
theorem. We observe first that Lemma 1.3 yields the following
special case of the theorem
$$
A=\rho_1 B^n \ , \quad
B=\rho_2 B^n \ , \quad
\Th=\{ (x,y)\in A\times B: x+y\in RB^n\}   \tag 1.2
$$
where $\rho_1,\rho_2,R > 0$ are arbitrary constants. The case
$\rho_1=1,\rho_2=\rho < 1$ follows directly and the general one by
symmetry and homogeneity.

\bigskip
The strategy for the rest of the proof is now as follows:
if $A_0,B_0\subset\Bbb R^n$ and $\Th_0\subset A_0\times B_0$, we
will show that there are $A,B,\Th$ of the form (1.2) verifying
\roster
\item"(i)" \ $\l(A_0)=\l(A), \ \ \l(B_0)=\l(B)$
\item"(ii)" \ $\l(\Th)\geq\l(\Th_0)$
\item"(iii)" $\l(A+_\Th B)\leq\l(A_0+_{\Th_0} B_0)$
\endroster
Now if the original $A_0,B_0,\Th_0$ had yielded a counterexample
to the theorem, the corresponding $A,B,\Th$ would have,
{\it a fortiori}, worked as such, contrary to the remark following
Lemma 1.3. Accordingly it remains to realize (i)--(iii) for given
$A_0,B_0,\Th_0$.

\medskip$\underline{\text{Step 1}}$. Set $C=A_0+_{\Th_0}B_0$
and
$$
\Th_1=\{(x,y)\in A_0\times B_0: x+y\in C\} \ ,
$$
then $A_0+_{\Th_0}B_0=A_0+_{\Th_1}B_0$, while clearly
$\l(\Th_1)\geq\l(\Th_0)$.

\medskip$\underline{\text{Step 2}}$. Define
$\rho_1,\rho_2,R > 0$ via
$$
\l(A_0)=\l(\rho_1 B^n) \ , \quad
\l(B_0)=\l(\rho_2 B^n) \ , \quad
\l(C)=\l(RB^n) \ .
$$
We then have
$$
\aligned
\l(\Th_1) & = \l(\{(x,y)\in A_0\times B_0: x+y\in C\}) \\
  & =\int_{\Bbb R^n}\int_{\Bbb R^n} \chi_{A_0}(x)\chi_{B_0}(y)
     \chi_C(x+y)dx dy \\
  & \leq \int_{\Bbb R^n}\int_{\Bbb R^n}\chi_{\rho_1 B^n}(x)
    \chi_{\rho_2 B^n}(y) \chi_{RB^n}(x+y)dxdy \\
  & = \l(\{(x,y)\in \rho_1B^n\times\rho_2B^n: x+y\in RB^n\})
\endaligned     \tag 1.3
$$
as required for (i)--(iii) (and concluding the derivation of
Theorem 1.2 from Lemma 1.3). The inequality in (1.3) is a special
case of [2, Theorem 3.4], which, in a much more general setting,
estimates an integral of a product of nonnegative functions by that
of their spherical (or Schwartz) symmetrizations;  we thank Alain Pajor
for pointing the paper [2] to us.\qed

\bigskip
{\smc Proof of Lemma 1.3} (Sketch). We will show that, for an
appropriate choice of $c_1 > 0$ and with
$\tau=\frac 12\min\{\rho\sqrt{n},1\}$, one has
$$
1\geq |x_0|\geq 1-\tau/n\Rightarrow
\l (\{ y: |y|\leq\rho, \ |x_0+y| >(1+\rho^2)^{\frac 12}\}) \geq
  c_1\l(\rho B^n)   \tag 1.4
$$
It then follows that
$$
\l(B^n\times\rho B^n\backslash\Th)\geq
(1-\tau/n)^n c_1\l(B^n)\cdot\l(\rho B^n)
$$
and that clearly implies the lemma. To show (1.4), we denote
$r_0=|x_0|$ and assume, as we may, that
$x_0=(r_0,0,\dots ,0)$ and $n\geq 2$. Then (the reader is advised
to draw a picture)
$$
\multline
|\{y: |y|\leq\rho, \ |x_0+y|\leq (1+\rho^2)^{\frac 12}\}| \\
  = |B^{n-1}|\cdot \left( \int^s_{-\rho}
  (\rho^2-u^2)^{\frac{n-1}{2}} du+\int^t_s
  (1+\rho^2-(r_0+u)^2)^{\frac{n-1}{2}}du\right)
\endmultline
$$
where $s=(1-r^2_0)/2r_0$ and $t=(1+\rho^2)^{\frac 12}-r_0$.
Since $s\leq (\tau/n)\cdot (1+\frac{\tau/n}{2r_0})\leq
(\rho/\sqrt{n})(1+O(n^{-1}))$, the contribution of the first
integral constitutes a proportion of $\l(\rho B^n)$ that is
strictly smaller than 1 (uniformly in $n$) and asymptotically, as
$n\to\infty$, is of order $\Phi(1)\cdot\l (\rho B^n)$,
where $\Phi$ is the c.d.f. of a standard $N(0,1)$ Gaussian random
variable. Similarly, the contribution of the second integral is
shown to be $o(1)\cdot\l(\rho B^n)$ as $n\to\infty$ (or, more
exactly, less than $(\rho/\sqrt{n})\cdot\l(\rho B^n)$ for all
$n\geq 2$); we omit the rather routine details. Combining the two
estimates yields (1.4), hence Lemma 1.3.

\bigskip
{\bf 1.4  Remark.} Theorem 1.2 is optimal in the following sense:
there exist constants $\a,A > 0$ such that, for any $n\in\Bbb N$
(resp. for any $n\in\Bbb N$, $\rho\in (0,1)$), there exist
$A,B\in\Bbb R^n$ (resp. with
$\rho\leq(\l(B)/\l(A))^{1/n}\leq\rho^{-1}$) and
$\Th\subset A\times B$ with $\l(\Th) >\a\l(A)\l(B)$
(resp. $\l(\Th) > (1-A\rho n^{\frac 12})\l(a)\l(B)$) such that
the assertion of the theorem does not hold.

\proclaim{1.5 \ Corollary} There exist $c,C > 0$ such that,
for any $\d\in [0,c]$, $n\in\Bbb N$, any $A,B\subset\Bbb R^n$
and any $\Th\subset A\times B$ with $\l(\Th)\geq (1-\d)\l(A)\l(B)$
one has
$$
\l(A+_{\Th} B)^{2/n}\geq (1-\frac{C\d}{n})(\l(A)^{2/n}+
\l(B)^{2/n})     \tag 1.5
$$
\endproclaim

{\smc Proof}. We may assume that $\l(A)=1\geq\l(B)=\rho^n$. Let
$c > 0$ be one given by Theorem 1.2;  we may clearly assume that
$c \leq 1/2$. If $\rho \geq \d/(c\sqrt{n})$, we may apply Thoerem 1.2
and get the assertion, in fact without the factor $(1-\frac{C\d}{n})$.
On the other hand,  regardless of the size of $\rho$ one has
(just by Fubini's theorem),
$$
\l(A +_{\Th} B)\geq (1-\d)\l(A) =1-\d,
$$
hence
$$
\l(A+_{\Th} B)^{2/n} \geq 1-\frac{3\d}{n} \ ,
$$
and it is easy to check that, for an appropriate choice of $C$,
the right-hand side of (1.5) does not exceed the latter quantity
if $\rho < \d/(c\sqrt{n})$.

\bigskip{\bf 1.6 \ Remark}. Redoing the argument of Theorem 2.1 in
the context of Corollary 1.5 (rather than formally applying the
assertion of the theorem) does not produce a sharper result.
However, it is possible to obtain an assertion similar to that of
Corollary 1.5 under much weaker assumptions, namely, in the
notation of Theorem 1.2, if $\g\in (0,1)$ then the condition
$\l(\Th)\geq\g\l(A)\l(B)$ implies a version of (1.5) with
$(1-C\d/n)$ replaced by  $(1-C\rho(\log(1+1/\g)/n)^{\frac 12})$.

\vskip .33 truein
\specialhead{2. The free entropy power inequality}
\endspecialhead

The free entropy $\chi(X_1,\dots ,X_n)$ for an $n$-tuple of
selfadjoint elements $X_j\in M$, $M$ a von Neumann alegbra with a
normal faithful trace state $\tau$, was defined in [8] part II.
The definitionn involves sets of matricial microstates
$\Gamma_R(X_1,\dots ,X_n;m,k,\e)$ (see $\S$2 in [8] part II).
The microstates are points in $(\Cal M^{sa}_k)^n$, where
$\Cal M^{sa}_k$ denotes the selfadjoint $k\times k$ matrices.
$\l$ will denote Lebesgue measure on $(\Cal M^{sa}_k)^n$
corresponding to the euclidean norm
$$
\|(A_1,\dots ,A_n)\|^2_{HS} = {\op{Tr}}(A^2_1 +\dots +A^2_n) \ .
$$
For one random variable we have (Prop. 4.5 in [8] part II) that:
$$
\chi(X)=\iint \log |s-t|d\mu(s)d\mu(t)+\tfrac 34 +\tfrac 12\log 2\pi
    \tag 2.1
$$
where $\mu$ is the distribution of $X$ (see 2.3 in [9]) or
equivalently the measure on $\Bbb R$ obtained by applying the trace
$\tau$ to the spectral measure of $X$.

\proclaim{2.1 \ Theorem} Let $X,Y\in M$, $X=X^*$, $Y=Y^*$
and assume $X,Y$ are free. Then
$$
\exp(2\chi(X))+\exp(2\chi(Y))\leq\exp(2\chi(X+Y)) \ .   \tag 2.2
$$
\endproclaim

Using the explicit formula for $\chi(X)$ and the fact that the
distribution of the sum of two free random variables is obtained
via the free convolution $\boxplus$ (see 3.1 in [9]) there is an
equivalent form of the preceding theorem.

\proclaim{2.1$'$ Theorem} Let $\a,\beta$ be compactly supported
probability measures on $\Bbb R$. Then
$$
\aligned
& \exp(2\iint\log |s-t|d\a(s)d\a(t)) +
\exp(2\iint\log |s-t|d\beta(s)d\beta(t)) \\
& \quad \leq
\exp(2\iint\log|s-t|d(\a\boxplus\beta)(s)d(\a\boxplus\beta) (t) \ .
\endaligned        \tag 2.3
$$
\endproclaim

{\smc Proof of Theorem 2.1} The proof will be technically similar
to sections 4 and 5 of [8] part II. Let $Z\in M$, $Z=Z^*$
distributed according to Lebesgue measure on [0,1] and let
$U_1,U_2$ be unitaries with Haar distributions in $(M,\tau)$
and assume $Z,U_1,U_2$ are $*$-free. Let further
$h_1,h_2: [0,1]\to\Bbb R$ be $C^1$-functions with
$h'_1(t) > 0$, $h'_2(t) > 0$ for all $t\in [0,1]$. Remark that it
suffices to prove the theorem in case $X=U_1h_1(Z)U^*_1$,
$Y=U_2h_2(Z)U^*_2$ (i.e., the distributions of $X$ and $Y$
are the push-forwards by $h_1$ and $h_2$ of Lebesgue measure on
$[0,1]$). Indeed see $2^\circ$ in the proof of Proposition 4.5 in
[8] part II) there are sequences $h_{j,n}$ of functions as above,
such that
$$
\align
& \lim_{n\to\infty} \chi(U_1 h_{1,n}(Z)U^*_1)=\chi(X) \\
& \lim_{n\to\infty} \chi(U_2 h_{2,n}(Z)U^*_2)=\chi(Y)
\endalign
$$
and $\|h_{j,n}\|_\infty < R$ for some fixed constant $R$. Then
$$
\|U_1h_{1,n}(Z)U^*_1+U_2 h_{2,n}(Z)U^*_2\|\leq 2R
$$
and $U_1 h_{1,n} (Z)U^*_1+U_2 h_{2,n}(Z)U^*_2$ converges in
distribution to $X+Y$ because of our freeness assumptions.
By 2.6 in [8] part II we have
$$
\limsup_{n\to\infty} \chi(U_1 h_{1,n}(Z)U^*_1+U_2 h_{2,n}(Z)U^*_2)
\leq \chi(X+Y)
$$
and hence it suffices to prove Theorem 2.1, in case
$X=U_1 h_1(Z)U^*_1$, \ $Y=U_2h_2(Z)U^*_2$.

Like in 5.3 of [8] part II, let
$$
\Omega(h_j;k) =
\{ A\in\Cal M^{sa}_k\mid h_j(2s/2k)\leq \l_{s+1}(A) \leq
h_j((2s+1)/2k), 0\leq s\leq k+1\}
$$
where $\l_1(A)\leq\dots\leq\l_k(A)$ are the eigenvalues of $A$.
The last part of the proof of Proposition 4.5 in [8] part II shows
that
$$
\lim_{k\to\infty} (k^{-2}\log\l(\Omega (h_j;k))+
2^{-1}\log k)=\chi(h_j(Z))    \tag 2.4
$$
where $\l$ is the Lebesgue measure on $\Cal M^{sa}_k$.

Let further $N\in\Bbb N$ and $\e > 0$ be given and
$$
\Th(k) = \{(A_1,A_2)\in\!\!\prod_{1\leq j\leq 2} \!\!
\Omega(h_j;k)\mid
(A_1,A_2)\in\Gamma (U_1 h_1(Z)U^*_1,U_2 h_2(Z)U^*_2;N,k,\e)\}
$$
By Lemma 5.3 in [8] part II we have:
$$
\lim_{k\to\infty} \
\frac{\l(\Th(k))}{\l(\Omega(h_1;k)\times\Omega(h_2;k))}=1  \tag 2.6
$$
If $R > \|h_j\|_{\infty}$ then
$$
\Th(k)\subset
\Gamma_R(U_1h_1(Z)U^*_1, U_2h_2(Z)U^*_2; \ N,k,\e) \ .
$$
Further, given $N_1\in\Bbb N$, \ $\e_1 > 0$ we may choose
$N\in\Bbb N$, \ $\e > 0$ so that
$$
(A_1,A_2)\in\Gamma_R(U_1h_1(Z)U^*_1,U_2h_2(Z)U^*_2; \ N,k,\e)
$$
implies
$$
A_1+A_2\in\Gamma_{2R}
 (U_1h_1(Z)U^*_1+U_2h_2(Z)U^*_2; \ N,k,\e) \ .
$$
In particular,
$$
\Omega(h_1;k) +_{\Th(k)} \Omega(h_2;k)
\subset \Gamma_{2R}(U_1h_1(Z)U^*_1+U_2h_2(Z)U^*_2; \ N_1,k,\e_1) \ .
   \tag 2.7
$$
Using Theorem 1.2 for $k\geq k_0$ with $k_0$ sufficiently large,
taking into account (2.6), we have
$$
(\l(\Omega(h_1;k)))^{2/k^2} + (\l(\Omega(h_2;k)))^{2/k^2}
\leq (\l(\Gamma_{2R}(U_1h_1(Z)U^*_1+U_2h_2(Z)U^*_2; \
N_1,k,\e_1)))^{2/k^2} \ .
$$
Given $\d > 0$ we may choose $k_0,N_1$ large and $\e_1$ small, so
that
$$
\align
& k^{-2}\log\l(\Gamma_{2R}(U_1h_1(Z)U^*_1+U_2h_2(Z)U^*_2; \
N_1,k,\e_1) + \tfrac 12\log k \\
& \qquad\leq \chi(U_1h_1(Z)U^*_1+U_2h_2(Z)U^*_2)+\d
\endalign
$$
for all $k\geq k_0$.

We infer that for $k\geq k_0$,
$$
\align
& \exp(2k^{-2}(\log\l(\Omega(h_1;k))+2^{-1}\log k)
 + \exp(2k^{-2}(\log\l(\Omega(h_2;k))+2^{-1}\log k) \\
& \qquad \leq \exp(2(\chi(U_1h_1(Z)U^*_1+U_2h_2(Z)U^*_2)+\d)) \ .
\endalign
$$
Letting $k\to\infty$ and taking into account that $\d > 0$
was arbitrary, we get the desired inequality.   \qed

\vskip .33 truein
\specialhead{3. Concluding remarks and open problems}\endspecialhead

{\bf 3.1 \ The free entropy power inequality for $n$-tuples}.
To extend Theorem 2.1 to $n$-tuples of non-commutative random
variables means to prove
$$
\exp(\frac 2n\chi(X_1,\dots ,X_n))+
\exp(\frac 2n\chi(Y_1,\dots ,Y_n))\leq
\exp(\frac 2n\chi(X_1+Y_1,\dots ,X_n+Y_n))  \tag 3.1
$$
under the assumption that $\{X_1,\dots ,X_n\}$ and
$\{Y_1,\dots ,Y_n\}$ are freee. The missing ingredient at this time
is the generalization of section 5 in [8] part II to $n$-tuples.
The rest of the argument, i.e. the use of Theorem 1.2, would then
be along the same lines as for $n=1$. At present, partial
generalizations of Theorem 2.1 can be obtained. The route to be
followed is: first replace $X$ and $Y$ by $n$-tuples
$(X_1,\dots ,X_n)$, $(Y_1,\dots ,Y_n)$ such that the $2n$ variables
$X_1,\dots ,X_n, \ Y_1,\dots ,Y_n$ are free and note that in this
situation the necessary facts about sets of matricial microstates
can be obtained from section 5 of [8] part II. Then the
generalization of Theorem 2.1 will hold for $n$-tuples
$(F_1(X_1,\dots,X_n),\dots ,F_n(X_1,\dots ,X_n))$ and
$(H_1(Y_1,\dots ,Y_n),\dots ,H_n(Y_1,\dots ,Y_n))$ where
$X_1,\dots ,X_n$, $Y_1,\dots ,Y_n$ are free and\newline
$(F_1,\dots ,F_n)$,$(H_1,\dots ,H_n)$ are non-commutative functions
satisfying suitable conditions, like the existence of an inverse of
the same kind and extending to the matricial microstates.
These kind of extensions have statements containing many technical
conditions, the proof, except for some technicalities, being along
the same lines as for $n=1$. We don't pursue this here, hoping that
better techniques will yield a proof of the free entropy power
inequality in full generality.

\bigskip
{\bf 3.2 \ Shannon's classical entropy power inequality and
restricted Minkowski sums}. We would like to signal that the
inequality in Theorem 1.2 has the potential to provide a proof also
of Shannon's classical entropy power inequality.  The reason is
that the classical entropy of an $n$-tuple of commutative random
variables can be defined via microstates (using the diagonal
subalgebra of $n\times n$ matrix algebra instead of the full
algebra) and the entropy power inequality would then
correspond to the same kind fo geometric problem at the level of
microstates as in the free case. We are thinking of exploring this
possibility in future work.

\bigskip
{\bf 3.3 \ The free analogue of the Stam inequality}. It seems
natural to look also for a free analogue of the Stam inequality
([7], see also [1],[3]), of which the free entropy power inequality
would be a consequence. With $\Phi$ denoting the free analogue of
Fisher's information measure (see [8] part I) this would amount to:
$$
(\Phi(X+Y))^{-1} \leq (\Phi(X))^{-1}+(\Phi(Y))^{-1}
$$
if $X,Y$ are free.

\newpage
\centerline{\bf References}
\roster
\item"1." Blachman, N. M.: The convolution inequality for entropy
  powers, {\it IEEE Trans. Inform. Theory} {\bf 2} (1965), 267--327.
\item"2." Brascamp, H. J., Lieb, E. H. and Luttinger, J. M.:
  A general rearrangment inequality for multiple integrals,
  {\it J. Funct. Analysis} {\bf 17} (1974), 227--237.
\item"3." Carlen, E. A.: Superadditivity of Fisher's information
  and logarithmic Sobolev inequalities, {\it J. Funct. Analysis}
  {\bf 101} (1991), 194--211.
\item"4." Pisier, G.: {\it The Volume of Convex Bodies and Banach
  Space Geometry}, Cambridge Univ. Press, 1989.
\item"5." Schneider, R.: Convex bodies: the Brunn-Minkowski
  theory, {\it Encyclopedia of Mathematics and Its Applications}
  {\bf 44}, Cambridge Univ. Press, 1993.
\item"6." Shannon, C. E. and Weaver, W.: {\it The Mathematical
  Theory of Communications}, University of Illinois Press, 1963.
\item"7." Stam, A. J.: Some inequalities satisfied by the
quantities of information of Fisher and Shannon, {\it Information
and Control} {\bf 2} (1959), 101--112.
\item"8." Voiculescu, D.: The analogues of entropy and of Fisher's
information measure in free probability theory, I, {\it Commun.
  Math. Phys.} {\bf 155} (1993), 71--92; \ {\it ibidem} II,
  {\it Invent. Math.} {\bf 118} (1994), 411--440; \
  {\it ibidem} III, IHES preprint, May 1995 (to appear in
  {\it GAFA}).
\item"9." Voiculescu, D., Dykema, D. and Nica, A.:
  {\it Free Random Variables}, CRM Monograph Series, vol. 1,
  American Mathematical Society, 1992.
\endroster
\enddocument